\documentclass{tbimc}

\newtheorem{thr}{Theorem}[section]
\newtheorem{lmm}{Lemma}[section]
\newtheorem{ass}{Assertion}[section]
\theoremstyle{definition}
\newtheorem{dfn}{Definition}[section]

\newcommand{\bydef}{\stackrel{\text{\rm def}}{=}}
\renewcommand{\P}{\mathcal P}
\DeclareMathOperator{\bcdot}{\boldsymbol\cdot}
\DeclareMathOperator{\re}{\rm Re}

\begin{document}

\selectlanguage{english}

\title[Continuous counterparts of Poisson and binomial distributions]
	{Continuous counterparts\\of Poisson and binomial distributions\\and their properties}

\author{Andrii Ilienko}
\address{Department of analysis and probability,
	National technical university of Ukraine \lq\lq Kiev polytechnic institute\rq\rq,
	37 Peremohy ave,
	Kiev, Ukraine}
\email{an.ilienko@gmail.com}


\begin{abstract}
On the basis of integral representations of Poisson and binomial distribution functions via complete and incomplete Euler $\Gamma$- and $\rm B$-functions, we introduce and discuss continuous counterparts of the Poisson and binomial distributions. The former turns out to be closely related to classical Volterra functions as well. Under usual condition $Np\to\lambda$, we also prove that the sequence of continuous binomial distributions converges weakly to the continuous Poisson one. At the end, we discuss a relationship between the continuous Poisson distribution and the $\Gamma$-process.
\end{abstract}

\maketitle

\section{Introduction}

In various applied research papers, many authors extensively use what they call a \lq\lq continuous Poisson distribution\rq\rq\ and a \lq\lq continuous binomial distribution\rq\rq, providing these terms with very different, not always correct meanings. For example, by the \lq\lq continuous Poisson distribution\rq\rq, some authors (\cite{Kim}, \cite{Tur}) understand an absolutely continuous distribution with the density of the form
$$f_\lambda(x)=C_\lambda\frac{e^{-\lambda}\lambda^x}{\Gamma(x+1)},\qquad x\ge0,$$
where $C_\lambda$ is a normalizing constant. Some others use to this end the $\Gamma$-distribution (\cite{Alis}, \cite{Herz}) or even simply the exponential distribution (\cite{Webb}).

From the strictly formal point of view, the above distributions can not be regarded as genuine continuous analogues of the classical Poisson law, since these distributions have probabilistically little in common with this law. In this paper, we make an attempt to introduce continuous counterparts of Poisson and binomial distributions in a more natural way, on the basis of integral representations of Poisson and binomial distribution functions via complete and incomplete Euler $\Gamma$- and $\rm B$-functions, and examine some properties of these counterparts.

Note that these distributions appeared in various forms (although without any detailed exposition) in the papers of G. Marsaglia \cite{Mars}, H. Stern \cite{Stern}, as well as in the solution of problem 3.4.1.22 in the book of D. Knuth \cite{Knuth}. In paper \cite{Mars} and book \cite{Knuth} they were employed for the computer simulation of the classical Poisson and binomial distributions. In report \cite{Stern}, it is noted that the continuous Poisson distribution is closely related to the distribution of the time when $\Gamma$-process jumps over a fixed level. This relation will be discussed in detail at the end of the paper.

In the present work, we proceed with the research started in paper \cite{Il}.

\section{Notation and preliminaries}

Recall the standard notation which will be used in what follows: for $x\in\mathbb R$ set
$$\lfloor x\rfloor\bydef\max\{k\in\mathbb Z\colon k\le x\},\quad
\lceil x\rceil\bydef\min\{k\in\mathbb Z\colon k\ge x\}.$$
Besides, we will need incomplete Euler $\Gamma$- and $\rm B$-functions:
\begin{align*}
\Gamma(x,\lambda)&\bydef\int_\lambda^\infty e^{-t}t^{x-1}\,dt, &x&>0,\,\lambda\ge0,
&&\text{(incomplete $\Gamma$-function)}\\
{\rm B}(x,y,p)&\bydef\int_p^1 t^{x-1}(1-t)^{y-1}\,dt, &x,y&>0,\,0\le p\le 1.
&&\text{(incomplete $\rm B$-function)}
\end{align*}
In particular, we have for the usual (complete) $\Gamma$- and $\rm B$-functions:
$$\Gamma(x)=\Gamma(x,0), \quad {\rm B}(x,y)={\rm B}(x,y,0), \quad x,y>0.$$

Denote by $\pi_\lambda$, $\lambda>0$, the classical Poisson distribution, i.e. a discrete probability measure on the set $\mathbb N\cup\{0\}$ given by $\pi_\lambda\bigl(\{k\}\bigr)=e^{-\lambda}\frac{\lambda^k}{k!}$. Also denote by $\beta_{N,p}$, $N\in\mathbb N$, $p\in[0,1]$, the classical binomial distribution, i.e. a discrete probability measure on the set $\{0,\dots,N\}$ given by $\beta_{N,p}\bigl(\{k\}\bigr)=\binom Nk p^k(1-p)^{N-k}$.

It is well-known that the distribution function $F_\lambda$ of the Poisson measure $\pi_\lambda$ admits a representation in terms of functions $\Gamma(\bcdot,\lambda)$ and $\Gamma(\bcdot)$:
\begin{equation}
\label{df_dp}
F_\lambda(x)\bydef\pi_\lambda\bigl((-\infty,x)\bigr)=
\begin{cases}
0,&x\le 0,\\e^{-\lambda}\sum_{k=0}^{\lceil x\rceil-1}\frac{\lambda^k}{k!},&x>0,
\end{cases}
=
\begin{cases}
0,&x\le 0,\\\frac{\Gamma\bigl(\lceil x\rceil,\lambda\bigr)}{\Gamma\bigl(\lceil x\rceil\bigr)},&x>0.
\end{cases}
\end{equation}
A similar representation via $\rm B$-functions holds true for the distribution function $F_{N,p}$ of the binomial measure $\beta_{N,p}$:
\begin{equation}
\begin{aligned}
\label{df_db}
F_{N,p}(x)\bydef\beta_{N,p}\bigl((-\infty,x)\bigr)&=
\begin{cases}
0,&x\le 0,\\\sum_{k=0}^{\lceil x\rceil-1}\binom Nk p^k(1-p)^{N-k},&0<x\le N,\\1,&x>N,
\end{cases}
\\&=
\begin{cases}
0,&x\le 0,\\\frac{{\rm B}\bigl(\lceil x\rceil,N+1-\lceil x\rceil,p\bigr)}
{{\rm B}\bigl(\lceil x\rceil,N+1-\lceil x\rceil\bigr)},&0<x\le N,\\1,&x>N.
\end{cases}
\end{aligned}
\end{equation}

Representations (\ref{df_dp}) and (\ref{df_db}) can be attributed to a sort of \lq\lq probabilistic folklore\rq\rq\ (see e.g. \cite{Shir}, where they are given as problems). The proofs are a simple exercise in integration by parts.

\section{Continuous analogues of Poisson and binomial distributions}

Taking representations (\ref{df_dp}) and (\ref{df_db}) into account, the following definitions appear to be quite natural.

\begin{dfn}
\label{cp}
By \emph{continuous Poisson distribution with parameter $\lambda>0$} we will mean the probability measure $\tilde\pi_\lambda$ supported by $[0,\infty)$ with distribution function of the form
\begin{equation}
\label{df_cp}
\tilde F_\lambda(x)=
\begin{cases}
0,&x\le 0,\\\frac{\Gamma(x,\lambda)}{\Gamma(x)},&x>0.
\end{cases}
\end{equation}
\end{dfn}

\begin{dfn}
\label{cb}
By \emph{continuous binomial distribution with parameters $N>0$, $p\in(0,1)$} we will mean the probability measure $\tilde\beta_{N,p}$ supported by $[0,N+1]$ with distribution function of the form
\begin{equation}
\label{df_cb}
\tilde F_{N,p}(x)=
\begin{cases}
0,&x\le 0,\\\frac{{\rm B}(x,N+1-x,p)}
{{\rm B}(x,N+1-x)},&0<x\le N+1,\\1,&x>N+1.
\end{cases}
\end{equation}
\end{dfn}

Note that in Definition \ref{cb} parameter $N$ may take non-integer values; we continue to denote it with $N$ only by tradition.

We should first check whether the distributions introduced by Definitions \ref{cp} and \ref{cb} are well defined.
To this end, we show that the following assertions hold true:
\begin{itemize}
\item[(i)] $\lim_{x\to +0}\frac{\Gamma(x,\lambda)}{\Gamma(x)}=0$;
\item[(ii)] $\lim_{x\to+\infty}\frac{\Gamma(x,\lambda)}{\Gamma(x)}=1$;
\item[(iii)] the function $\frac{\Gamma(x,\lambda)}{\Gamma(x)}$ increases in $x$ on the interval $[0,\infty)$;
\item[(iv)] $\lim_{x\to +0}\frac{{\rm B}(x,N+1-x,p)}{{\rm B}(x,N+1-x)}=0$;
\item[(v)] $\lim_{x\to N+1-0}\frac{{\rm B}(x,N+1-x,p)}{{\rm B}(x,N+1-x)}=1$;
\item[(vi)] the function $\frac{{\rm B}(x,N+1-x,p)}{{\rm B}(x,N+1-x)}$ increases in $x$ on the segment $[0,N+1]$.
\end{itemize}

Assertions (i) and (iv) hold, since the numerators are bounded as $x\to+0$ while the denominators increase to infinity.  Similarly, ${\rm B}(x,N+1-x)-{\rm B}(x,N+1-x,p)$ is bounded as $x\to N+1-0$ while ${\rm B}(x,N+1-x)$ increases to infinity, which yields (v). (ii) follows from the fact that $\Gamma(x)-\Gamma(x,\lambda)$ has an asymptotic order $x^{-1}\lambda^x$ as $x\to+\infty$, in other words is $o\bigl(\Gamma(x)\bigr)$.

Finally, assertions (iii) and (vi) result from the following relations:
\begin{gather}
\begin{aligned}
\label{d_cp}
\tilde f_\lambda(x)\bydef \tilde F'_\lambda(x)
&=\Gamma^{-2}(x)\biggl(\int_\lambda^\infty\int_0^\infty e^{-(s+t)}(st)^{x-1}\ln\frac st\,ds\,dt\biggr)
\\&=\Gamma^{-2}(x)\biggl(\int_\lambda^\infty\int_0^\lambda e^{-(s+t)}(st)^{x-1}\ln\frac st\,ds\,dt\biggr)>0,\qquad x>0,
\end{aligned}\\
\begin{aligned}
\label{d_cb}
\tilde f_{N,p}(x)&\bydef \tilde F'_{N,p}(x)
\\={\rm B}^{-2}&(x,N+1-x)\biggl(\int_p^1\int_0^1 (st)^{x-1}(1-s)^{N-x}(1-t)^{N-x}\ln\frac{s(1-t)}{t(1-s)}\,ds\,dt\biggr)
\\={\rm B}^{-2}&(x,N+1-x)\biggl(\int_p^1\int_0^p (st)^{x-1}(1-s)^{N-x}(1-t)^{N-x}\ln\frac{s(1-t)}{t(1-s)}\,ds\,dt\biggr)>0,
\\&\hspace{9.1 cm}{0<x<N+1.}
\end{aligned}
\end{gather}
Note that we changed $\infty$ for $\lambda$ and $1$ for $p$ in the upper limits of the inner integrals due to antisymmetry of the integrands with respect to the argument pair $(s,t)$.

The result can be summarized in the next statement. Here $\mathbb R^+$ stands for $[0,\infty)$.

\begin{thr}\label{d_cpb}
The distribution functions $\tilde F_\lambda$ and $\tilde F_{N,p}$ given by (\ref{df_cp}) and (\ref{df_cb}) define absolutely continuous (with respect to the Lebesgue measure) probability distributions on the measurable space
$\bigl(\mathbb R^+,\mathfrak B\bigl(\mathbb R^+\bigr)\bigr)$. Their densities are of form (\ref{d_cp}) and (\ref{d_cb}).
\end{thr}
\pagebreak

\section{Volterra functions and the moments\\of the continuous Poisson distribution}

In this section, we study the moments of the continuous Poisson distribution. They turn out to be closely related to classical Volterra functions, which play a role in integral transforms.

We first give the corresponding definitions (see e.g. \cite{Apel1} or \cite{Apel2}). In the most general form, the Volterra $\mu$-function is defined by
$$\mu(t,\alpha,\beta)=\int_0^\infty\frac{t^{x+\alpha}x^\beta}{\Gamma(\beta+1)\Gamma(x+\alpha+1)}\,dx,\qquad \re\beta>-1.$$
In particular, by setting $\alpha=\beta=0$ we come to a simpler and more known Volterra $\nu$-function: $\nu(t)=\int_0^\infty\frac{t^x}{\Gamma(x+1)}\,dx$.

For a fixed $k\in\mathbb N$, consider the $k$-order moment function of the continuous Poisson distribution:
$$m_k(\lambda)\bydef\int_0^\infty x^k\tilde\pi_\lambda(dx),\qquad\lambda>0.$$
Then we have
$$m_k(\lambda)=k\int_0^\infty x^{k-1}(1-\tilde F_\lambda(x))\,dx=k\int_0^\infty \frac{x^{k-1}\bigl(\Gamma(x)-\Gamma(x,\lambda)\bigr)}{\Gamma(x)}\,dx.$$
By the definition of complete and incomplete $\Gamma$-functions and Fubini theorem
\begin{equation}
\label{pm_v}
m_k(\lambda)=k\int_0^\lambda e^{-t}\biggl(\int_0^\infty\frac{x^{k-1}t^{x-1}}{\Gamma(x)}\,dx\biggr)\,dt=
k!\int_0^\lambda e^{-t}\mu(t,-1,k-1)\,dt.
\end{equation}

The moments of the continuous Poisson distribution look particularly nice in terms of Laplace transform in $\lambda$. Again by Fubini theorem we have
\begin{equation}
\label{pm_l1}
\hat m_k(s)\bydef\int_0^\infty e^{-s\lambda}m_k(\lambda)\,d\lambda=\frac{k!}s\int_0^\infty e^{-t(s+1)}\mu(t,-1,k-1)\,dt.
\end{equation}
The Laplace transform $\hat\mu$ of the Volterra function $\mu$ is of the form
$$\hat\mu(s,\alpha,\beta)=\frac 1{s^{\alpha+1}(\ln s)^{\beta+1}},\qquad\re s>1.$$
(See e.g. \cite{Apel2}; anyway, this can be deduced by a simple calculation.) Thus, (\ref{pm_l1}) yields
\begin{equation}
\label{pm_l2}
\hat m_k(s)=\frac{k!}{s\ln^k(1+s)},\qquad\re s>0.
\end{equation}

The previous relation makes it possible to find the double Laplace transform (Laplace-Stieltjes with respect to measure and Laplace with respect to $\lambda$) of the distribution family $(\tilde\pi_\lambda,\lambda>0)$. Let $\xi_\lambda$, $\lambda>0$, denote a random variable distributed by the continuous Poisson law with parameter $\lambda$. Then, for $\re s>0$ we have
\begin{equation}
\label{dlt}
\begin{aligned}
&\hat\varphi(p,s)\bydef\int_0^\infty e^{-s\lambda}\mathbb Ee^{-p\xi_\lambda}\,d\lambda=
\iint_{(0,\infty)\times(0,\infty)}e^{-px-s\lambda}\,\tilde\pi_\lambda(dx)\,d\lambda=\\
&=\sum_{k=0}^\infty\frac{(-p)^k}{k!}\int_0^\infty e^{-s\lambda}m_k(\lambda)\,d\lambda=
\sum_{k=0}^\infty\frac{(-p)^k\hat m_k(s)}{k!}=\frac 1s\frac{\ln(1+s)}{p+\ln(1+s)}.
\end{aligned}
\end{equation}

We now resume all the above in the next statement.

\begin{thr}
The moments of the continuous Poisson distribution $\tilde\pi_\lambda$ are given by formula (\ref{pm_v}). In terms of Laplace transform they can be expressed as (\ref{pm_l2}). The double Laplace transform of the distribution family $(\tilde\pi_\lambda,\lambda>0)$ is of form (\ref{dlt}).
\end{thr}

\section {Convergence of continuous binomial distributions}

The classical Poisson theorem (sometimes also called the law of rare events) asserts that under the condition $Np\to\lambda$ the sequence of binomial distributions with parameters $N,p$ weakly converges to the Poisson one with parameter $\lambda$. In this section, we show that this remains true also for their continuous counterparts.

\begin{thr}
\label{ct}
Let $\bigl(\tilde\beta_{N(k),p(k)},k\in\mathbb N\bigr)$ be a sequence of continuous binomial distributions. If $N(k)\to\infty$ and $p(k)\to0$ as $k\to\infty$ in such a way that $N(k)p(k)\to\lambda$ then $\tilde\beta_{N(k),p(k)}\xrightarrow{w}\tilde\pi_\lambda$.
\end{thr}

We break the proof up into three lemmas.

\begin{lmm}
\label{l1}
For arbitrary $x\in[0,N]$ and $y\in[0,\infty)$
\begin{align}
\label{l11}
\tilde\beta_{N,p}\bigl([x,x+1)\bigr)&=\frac{\Gamma(N+1)}{\Gamma(x+1)\Gamma(N-x+1)}\cdot p^x(1-p)^{N-x},\\
\tilde\pi_\lambda\bigl([y,y+1)\bigr)&=\frac{e^{-\lambda}\lambda^y}{\Gamma(y+1)}.
\label{l12}
\end{align}
\end{lmm}

\emph{Proof.} Performing a simple integration by parts, we have
\begin{align*}
\tilde\beta_{N,p}\bigl([x,x+1)\bigr)=\tilde F_{N,p}(x+1)-\tilde F_{N,p}(x)&=
\frac{{\rm B}(x+1,N-x,p)}{{\rm B}(x+1,N-x)}-\frac{{\rm B}(x,N+1-x,p)}{{\rm B}(x,N+1-x)}=\\
&=\frac{p^x(1-p)^{N-x}}{(N-x){\rm B}(x+1,N-x)},
\end{align*}
which yields (\ref{l11}) by the relationship between $\Gamma$- and ${\rm B}$-functions. (\ref{l12}) can be obtained in a similar manner.

\begin{lmm}
\label{l2}
Under the conditions of Theorem \ref{ct}
$$\lim_{k\to\infty}\tilde\beta_{N(k),p(k)}\bigl([x,x+1)\bigr)=\tilde\pi_\lambda\bigl([x,x+1)\bigr),\qquad x\ge0.$$
\end{lmm}

The proof follows from Lemma \ref{l1} and elementary computations with Stirling's formula and is thus omitted.

\begin{lmm}
\label{l3}
Let $\bigl(\P_k,k\in\mathbb N\cup\{\infty\}\bigr)$ be a sequence of probability measures on the measurable space
$\bigl(\mathbb R^+,\mathfrak B\bigl(\mathbb R^+\bigr)\bigr)$. If for each $x\ge0$
\begin{equation}
\label{l2c}
\lim_{k\to\infty}\P_k\bigl([x,x+1)\bigr)=\P_\infty\bigl([x,x+1)\bigr),
\end{equation}
then the weak convergence of $\P_k$ to $\P_\infty$ holds.
\end{lmm}

Lemma \ref{l3} actually implies that the intervals $[x,x+1)$, $x\ge0$, form a convergence determining class (see e.g. \cite{Bil}, p. 18).

\medskip
\emph{Proof.} Denote by $F_k$ the distribution function of the measure $\P_k$. For a fixed $\varepsilon>0$ let $N_\varepsilon$ be such an integer that
$\P_\infty\bigl([N_\varepsilon,\infty)\bigr)<\varepsilon$. Then, for each $k\in\mathbb N\cup\{\infty\}$ and $x\ge0$
$$F_k(x)=1-\P_k\bigl([x,\infty)\bigr)=
1-\sum_{n=1}^{\lceil N_\varepsilon-x\rceil}\P_k\bigl([x+n-1,x+n)\bigr)
-\P_k\bigl([x+\lceil N_\varepsilon-x\rceil,\infty)\bigr),$$
and so
\begin{equation}
\begin{aligned}
\label{ddf}
|F_\infty(x)-&F_k(x)|\le\\\le
&\Biggl|\sum_{n=1}^{\lceil N_\varepsilon-x\rceil}\P_\infty\bigl([x+n-1,x+n)\bigr)-
\sum_{n=1}^{\lceil N_\varepsilon-x\rceil}\P_k\bigl([x+n-1,x+n)\bigr)\Biggr|+\\
+&\P_\infty\bigl([x+\lceil N_\varepsilon-x\rceil,\infty)\bigr)+
\P_k\bigl([x+\lceil N_\varepsilon-x\rceil,\infty)\bigr).
\end{aligned}
\end{equation}

The absolute value on the right-hand side in (\ref{ddf}) tends to zero as $k\to\infty$ due to condition (\ref{l2c}).
Since $x+\lceil N_\varepsilon-x\rceil\ge N_\varepsilon$, we have
$\P_\infty\bigl([x+\lceil N_\varepsilon-x\rceil,\infty)\bigr)<\varepsilon$. Moreover,
\begin{align*}
\limsup_{k\to\infty}\P_k\bigl([x+\lceil N_\varepsilon-x\rceil,\infty)\bigr)&\le
\lim_{k\to\infty}\P_k\bigl([N_\varepsilon,\infty)\bigr)=1-\lim_{k\to\infty}\P_k\bigl([0,N_\varepsilon)\bigr)=\\
&=1-\P_\infty\bigl([0,N_\varepsilon)\bigr)=\P_\infty\bigl([N_\varepsilon,\infty)\bigr)<\varepsilon.
\end{align*}
Letting now at first $k\to\infty$ and then $\varepsilon\to0$, we obtain by (\ref{ddf}) that
$$\lim_{k\to\infty}F_k(x)=F_\infty(x),\qquad x\ge0,$$
which was to be proved.

\medskip
Proof of Theorem \ref{ct} follows immediately from Lemmas \ref{l2} and \ref{l3}.

\section{An application to the $\Gamma$-process}

In this section, we show that the continuous Poisson distribution appears as the distribution of the (in a proper way normalized) time when $\Gamma$-process jumps over a fixed level.

We start with an informal explanation of this phenomenon. Consider a sequence of independent identically exponentially distributed random variables $\{\tau_i,i\in\mathbb N\}$ and the sequence of their cumulative sums $T=\{T_n,n\in\mathbb N\}$: $T_n=\sum_{i=1}^n\tau_i$. Also set $T_0=0$. It is well known that the renewal process $\bigl(X(t),t\ge0\bigr)$ constructed by the sequence $T$ as $X(t)=\sup\{n\colon T_n\le t\}$ is the classical Poisson process. The question is whether this setting can be carried into continuous time in such a way that the continuous Poisson distribution occurs instead of the classical one.

To answer this question we look upon the sequence $T$ as a L\'evy process in discrete time, i.e. random sequence with independent homogeneous increments. We can carry this process into continuous time, since the distribution of unit increments (the exponential one) is infinitely divisible. What's more, it is a special case of the $\Gamma$-distribution. Thus we may regard the $\Gamma$-process as a natural continuous counterpart of the sequence $T$.

Recall the corresponding definition (see e.g. \cite{Apple}, p. 54).

\begin{dfn}
A L\'evy process $\bigl(\tilde T(t),t\ge0\bigr)$ is called the $\Gamma$-process with parameters $\alpha,\beta>0$ if its transition density has a form
$$f_t(x)=\frac{\beta^{\alpha t}}{\Gamma(\alpha t)}x^{\alpha t-1}e^{-\beta x},\qquad x\ge0.$$
\end{dfn}

\begin{ass} Let $\tau_c$ be a time when the process $\tilde T$ jumps over a fixed level $c$. Then the random variable $\alpha\tau_c$ is continuously Poisson distributed with parameter $\beta c$.
\end{ass}

The proof follows from the fact that the distribution function $F_c$ of the random variable $\alpha\tau_c$ coincides with that of the continuous Poisson distribution $\tilde\pi_{\beta c}$:
\begin{align*}
F_c(x)&=\mathbb P\Bigl\{\tau_c<\frac x\alpha\Bigr\}=\mathbb P\Bigl\{\tilde T\Bigl(\frac x\alpha\Bigr)>c\Bigr\}
=\int_c^\infty f_{x/\alpha}(u)\,du=\frac{\beta^x}{\Gamma(x)}\int_c^\infty u^{x-1}e^{-\beta u}\,du=\\
&=\frac{\Gamma(x,\beta c)}{\Gamma(x)}=\tilde F_{\beta c}(x),\qquad x>0.
\end{align*}
\pagebreak

\end{document}